# Approximating the Analytic Fourier Transform with the Discrete Fourier Transform

Jeremy Axelrod[*]

*Department of Physics, University of California, Berkeley*

(Dated: 26 May 2015)

The Fourier transform is approximated over a finite domain using a Riemann sum. This Riemann sum is then expressed in terms of the discrete Fourier transform, which allows the sum to be computed with a fast Fourier transform algorithm more rapidly than via a direct matrix multiplication. Advantages and limitations of using this method to approximate the Fourier transform are discussed, and prototypical MATLAB codes implementing the method are presented.

### I. Introduction

The Fourier transform is a ubiquitous analytical mathematical tool. However, in many problems of interest analytic expressions for transformed functions do not exist, or the function to be transformed is only known at a set of discrete points as is the case for most real-world experimental data. In both of these cases, it is necessary to approximate the Fourier transform on a set of discrete points. This can be done by approximating the integral in the Fourier transform as a Riemann sum. Such a summation implemented as a single matrix multiplication by the vector of points to be transformed results in an undesirable algorithm complexity scaling of $\mathcal{O}(N^2)$. Here, an algorithm is presented which allows this Riemann sum to be expressed in terms of the discrete Fourier transform (DFT), which can in turn be computed via a fast Fourier transform (FFT) with complexity scaling of $\mathcal{O}(N \log N)$ [1]. The overall complexity of this algorithm also scales as $\mathcal{O}(N \log N)$.

### II. Method

#### A. Definitions

Let the Fourier transform, $\tilde{f}(\omega)$, of a function $f(t') : \mathbb{R} \to \mathbb{C}$ be defined as

$$\tilde{f}(\omega) := \sqrt{\frac{|b|}{(2\pi)^{1-a}}} \int_{-\infty}^{\infty} f(t') e^{ib\omega t'} dt' \quad (1)$$

where

$$f(t) := \sqrt{\frac{|b|}{(2\pi)^{1+a}}} \int_{-\infty}^{\infty} \tilde{f}(\omega') e^{-ib\omega' t} d\omega' \quad (2)$$

is the inverse Fourier transform of $\tilde{f}$, and $a, b \in \mathbb{R}$, $b \neq 0$ are arbitrary constants chosen by convention [2]. Let $\vec{x} := [x_1, x_2, \cdots, x_N]$ be a vector of length $N$, and let $j, j', k, k' = 1, 2, 3, ..., N$. The discrete Fourier transform, $\vec{X} := [X_1, X_2, \cdots, X_N]$, of $\vec{x}$ is then defined to be such that

$$X_k := \sum_{j'=1}^{N} x_{j'} e^{-\frac{2\pi i}{N}(j'-1)(k-1)} \quad (3)$$

where

$$x_j := \frac{1}{N} \sum_{k'=1}^{N} X_{k'} e^{\frac{2\pi i}{N}(j-1)(k'-1)}$$

is the inverse discrete Fourier transform of $\vec{X}$.

#### B. Approximating the forward transform

Let the function $f(t')$ be represented by a vector $\vec{x}$ of length $N$, so that $x_{j'} = f(t'_{j'})$ where $t'_{j'} = \tau'(j' - 1) + t'_1$ constitutes a vector of evenly-spaced points beginning at $t'_1$. Then (1) can be approximated via a Riemann sum, $\tilde{F}(\omega)$, as

$$\tilde{f}(\omega) \cong \tilde{F}(\omega) := \sqrt{\frac{|b|}{(2\pi)^{1-a}}} |\tau'| \sum_{j'=1}^{N} x_{j'} e^{ib\omega t'_{j'}} \quad (4)$$

Letting $\max\{t'_{j'}\}$, $\min\{t'_{j'}\}$ denote the maximum and minimum values in $\{t'_{j'}\}$, respectively, then

$$\lim_{N \to \infty} \tilde{F}(\omega) = \sqrt{\frac{|b|}{(2\pi)^{1-a}}} \int_{\min\{t'_{j'}\}}^{\max\{t'_{j'}\}} f(t') e^{ib\omega t'} dt'$$

when the limits of integration $\min\{t'_{j'}\}$ and $\max\{t'_{j'}\}$ are held fixed so that $|\tau'|$ decreases as $N$ increases. Therefore, $\tilde{F}(\omega)$ can be viewed as a Riemann approximation of the Fourier transform of the function

$$f_0(t') := \begin{cases} f(t') & \forall \, t' \in [\min\{t'_{j'}\}, \max\{t'_{j'}\}] \\ 0 & \text{elsewhere} \end{cases}$$

[*] jaxelrod@berkeley.edu

so that $\tilde{F}(\omega)$ approximates $\tilde{f}(\omega)$ well for large $N$, and for $\{t'_{j'}\}$ and $f(t')$ such that

$$\int_{-\infty}^{\min\{t'_{j'}\}} f(t')e^{ib\omega t'} dt' + \int_{\max\{t'_{j'}\}}^{\infty} f(t')e^{ib\omega t'} dt'$$
$$\ll \int_{\min\{t'_{j'}\}}^{\max\{t'_{j'}\}} f(t')e^{ib\omega t'} dt' \quad \forall \, \omega$$

Computing the sum in (4) directly for $N$ values of $\omega$ (denoted by $\vec{\omega} := [\omega_1, \omega_2, \cdots, \omega_N]$) can be written as a matrix multiplication:

$$\tilde{F}(\vec{\omega}) = \sqrt{\frac{|b|}{(2\pi)^{1-a}}} |\tau'| \begin{bmatrix} e^{ib\omega_1 t'_1} & e^{ib\omega_1 t'_2} & \cdots & e^{ib\omega_1 t'_N} \\ e^{ib\omega_2 t'_1} & e^{ib\omega_2 t'_2} & \cdots & e^{ib\omega_2 t'_N} \\ \vdots & \vdots & \ddots & \vdots \\ e^{ib\omega_N t'_1} & e^{ib\omega_N t'_2} & \cdots & e^{ib\omega_N t'_N} \end{bmatrix} \begin{bmatrix} x_1 \\ x_2 \\ \vdots \\ x_N \end{bmatrix} \quad (5)$$

Since multiplying a vector into an $N \times N$ matrix involves $N^2$ multiplication operations and $N(N-1)$ addition operations, the asymptotic complexity of computing (5) scales as $\mathcal{O}(N^2)$. However, "linearithmic" $\mathcal{O}(N \log N)$ scaling can be achieved by expressing $\tilde{F}(\vec{\omega})$ in terms of the DFT as defined in (3) and then utilizing an FFT algorithm to compute the DFT [1]. Referencing (3), in order to express $\tilde{F}(\omega)$ in terms of $\vec{X}$, define

$$\omega_k := -\frac{2\pi}{\tau' b N}(k-1) \quad (6)$$

so that

$$X_k = \sum_{j'=1}^{N} x_{j'} e^{ib\omega_k(t'_{j'} - t'_1)}$$
$$= e^{-ib\omega_k t'_1} \sqrt{\frac{(2\pi)^{1-a}}{|b|}} \frac{1}{|\tau'|} \tilde{F}(\omega_k)$$
$$\Rightarrow \tilde{F}(\omega_k) = \sqrt{\frac{|b|}{(2\pi)^{1-a}}} |\tau'| e^{ib\omega_k t'_1} X_k$$

As such, $\tilde{F}$ has been expressed in terms of the DFT, $\vec{X}$, for $\vec{\omega} = \{\omega_k\}$ as defined by (6). However, it is sometimes desirable to retrieve the values of $\tilde{F}$ for $\omega \notin \{\omega_k\}$, e.g. $\omega < 0$ if $b < 0$. To that end, note that for $m \in \mathbb{Z}$,

$$\tilde{F}\left(\omega + \frac{2\pi}{\tau' b} m\right) = \sqrt{\frac{|b|}{(2\pi)^{1-a}}} |\tau'| \sum_{j'=1}^{N} x_{j'} e^{ib\omega t'_{j'}} e^{i\frac{2\pi}{\tau'} m t'_{j'}}$$
$$= \sqrt{\frac{|b|}{(2\pi)^{1-a}}} |\tau'| \sum_{j'=1}^{N} x_{j'} e^{ib\omega t'_{j'}} e^{i2\pi m(j'-1)} e^{i\frac{2\pi}{\tau'} m t'_1}$$
$$= e^{i\frac{2\pi}{\tau'} m t'_1} \tilde{F}(\omega)$$

so that $\left|\tilde{F}\right|$ is periodic with period $\frac{2\pi}{\tau' b}$, with a phase shift of $\frac{2\pi}{\tau'} m t'_1$. Since $\vec{\omega}$ spans an entire period (from 0 to $\frac{-2\pi}{\tau' b} \frac{N-1}{N}$), $\tilde{F}$ can be determined from $\tilde{F}(\vec{\omega})$ for any integer multiple of $\frac{2\pi}{\tau' b N}$. This is useful in practice, since it is often desirable to retrieve $\tilde{F}$ on the interval $\omega \in [-|\omega_{nyq}|, +|\omega_{nyq}|]$ for $\omega_{nyq} := -\frac{\pi}{\tau' b}$ the Nyquist frequency in order to display aliases in a more readily interpretable context–that is, negative frequency components will appear below $\omega = 0$ instead of above $\omega = \omega_{nyq}$. A comparison between the DFT and the Riemann sum approximation to the Fourier transform is shown in figure 1. Clearly, if the vector being transformed is not entirely real, considering the DFT only below the Nyquist frequency index disregards possible asymmetries between positive and negative frequency components. Also, the vertical scaling of the Riemann sum approximation is independent of $N$ unlike with the DFT, making direct comparisons of spectral power density between transforms of vectors of different lengths possible.

### C. Approximating the inverse transform

The inverse transform can be treated in exactly the same way as the forward transform. Letting the function $\tilde{f}(\omega')$ be represented by a vector $\vec{X}$ of length $N$ so that $X_{k'} = \tilde{f}(\omega'_{k'})$ where $\omega'_{k'} = W'(k'-1) + \omega'_1$ constitutes a vector of evenly-spaced points beginning at $\omega'_1$, then (2) can be approximated via a Riemann sum, $F(t)$, as

$$f(t) \cong F(t) := \sqrt{\frac{|b|}{(2\pi)^{1+a}}} |W'| \sum_{k'=1}^{N} X_{k'} e^{-ib\omega'_{k'} t} \quad (7)$$

Then, after defining $t_j := -\frac{2\pi}{W'bN}(j-1)$,

$$F(t_j) = \sqrt{\frac{|b|}{(2\pi)^{1+a}}} N |W'| e^{-ib\omega'_1 t_j} x_j$$

and for $l \in \mathbb{Z}$,

$$F\left(t + \frac{2\pi}{W'b} l\right) = e^{-i\frac{2\pi}{W'} l \omega'_1} F(t)$$

similarly to before.

### D. Invertibility of the approximate transforms

It is of interest to know if the above Riemann sum approximations of the forward/inverse Fourier transforms are inverses of each other, i.e. is $F$ applied to $\vec{X} = \tilde{F}(\vec{\omega})$, where $\tilde{F}$ is applied to vector $\vec{x}$, equal to $\vec{x}$? Letting $x_{j'}$ be defined on $t'_{j'} = \tau'(j'-1) + t'_1$, and using equations (4) and (7), the inverse transform of the forward transform of $\vec{x}$ can be written as

$$F(t) = \sqrt{\frac{|b|}{(2\pi)^{1+a}}} |W'| \sum_{k'=1}^{N} \left( \sqrt{\frac{|b|}{(2\pi)^{1-a}}} |\tau'| \sum_{j'=1}^{N} x_{j'} e^{ib\omega'_{k'} t'_{j'}} \right) e^{-ib\omega'_{k'} t}$$

$$= \frac{|b|}{2\pi} |W'| |\tau'| \sum_{j'=1}^{N} x_{j'} e^{ib(t'_{j'}-t)\omega'_1} \left( \sum_{k'=1}^{N} \left( e^{ibW'(t'_{j'}-t)} \right)^{(k'-1)} \right) \quad (8)$$

assuming that $\omega'_{k'} = W'(k'-1) + \omega'_1$ in accordance with equation (6) and the accompanying discussion on discrete shifts. The sum in parentheses in (8) can be evaluated as a geometric sum:

$$\sum_{k'=1}^{N} \left( e^{ibW'(t'_{j'}-t)} \right)^{k'-1}$$
$$= \begin{cases} \frac{1-e^{iNbW'(t'_{j'}-t)}}{1-e^{ibW'(t'_{j'}-t)}} & \text{for } e^{ibW'(t'_{j'}-t)} \neq 1 \\ N & \text{for } e^{ibW'(t'_{j'}-t)} = 1 \end{cases}$$

In order to consider invertibility, it must be that $t \in \{t'_{j'}\}$. If this is the case, then $\forall\, j' \in \{1,2,3,...,N\}$, $t'_{j'} - t = m|\tau'|$ for some $m \in \{\mathbb{Z} \mid |m| \leq N-1\}$. That is, the difference between $t$ and $t'_{j'}$ is always an integer multiple of $|\tau'|$. Thus, since $|\tau'| = \frac{2\pi}{N|W'||b|}$,

$$e^{iNbW'(t'_{j'}-t)} = e^{\operatorname{sgn}(bW')i2\pi m} = 1$$
$$e^{ibW'(t'_{j'}-t)} = e^{\operatorname{sgn}(bW')i2\pi m/N} \neq 1$$
$$\forall\, m \in \{\mathbb{Z} \setminus 0 \mid |m| \leq N-1\}$$

Therefore,

$$\sum_{k'=1}^{N} \left( e^{ibW'(t'_{j'}-t)} \right)^{k'-1} = \begin{cases} 0 & \text{for } t'_{j'} \neq t \\ N & \text{for } t'_{j'} = t \end{cases}$$

and so

$$F(t) = \frac{|b|}{2\pi} |\tau'| |W'| N x_{j'} = x_{j'} \quad \text{when } t = t'_{j'}$$

which proves that the Riemann sum approximations are indeed inverses of each other. Since the method for calculating the Riemann sum approximation of the inverse Fourier transform using the inverse DFT can only return $t$ which are integer multiples of $\frac{2\pi}{NW'b}$, the transforms $F(t)$ and $\tilde{F}(\omega)$ specified above are restricted to being each others' inverse under the single condition that the input $\{t'_{j'}\}$ is of the form $t'_{j'} = \tau'(j' - 1 + l)$ so that it is possible for the inverse transform to return $t \in \{t'_{j'}\}$.

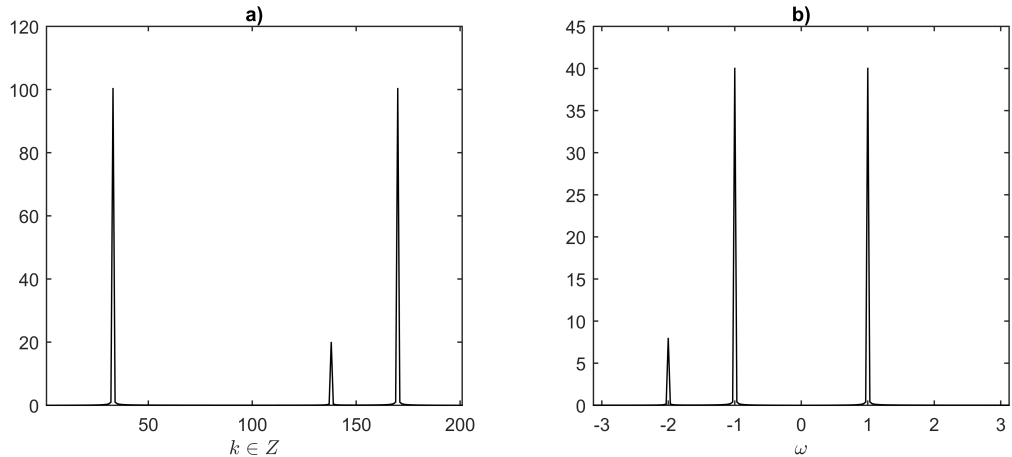

FIG. 1. Comparison of the magnitude of the DFT (a) and Riemann sum approximation to the Fourier transform (b) for $f(t) = \sin(t) + 0.1e^{-2it}$ sampled at 201 evenly-spaced points on the interval $[-100, 100]$ with $a = 0$, $b = -1$. The DFT is more difficult to directly interpret because of aliasing across the Nyquist frequency at $k = 101$.

### III. Performance

A prototypical MATLAB function, `FT`, which implements the approximation method described above for the forward Fourier transform is given in Section A1. The Riemann sum approximation to the Fourier transform of $f(t) = \text{rect}(t - 1)$ calculated using `FT` is shown in figure 2. The analytic Fourier transform of $f(t)$ for $a = 0$, $b = -1$ is

$$\tilde{f}(\omega) = \frac{1}{\sqrt{2\pi}} e^{-i\omega} \frac{\sin\left(\frac{\omega}{2}\right)}{\omega/2} \qquad (9)$$

and is also plotted in figure 2 for comparison. Referencing the code in `FT`, it is evident that adapting the DFT to the Riemann sum approximation only incurs additional computational costs scaling as $\mathcal{O}(N)$ from the element-wise multiplications by phase factors. Since the FFT complexity scales as $\mathcal{O}(N \log N)$, the overall complexity of the adaptation algorithm scales as $\mathcal{O}(N \log N)$. This is demonstrated in figure 3.

### IV. Conclusions

It has been shown that the Riemann sum approximation to the Fourier integral over a discrete finite domain can be expressed in terms of the discrete Fourier transform and is therefore calculable using a fast Fourier transform algorithm, which reduces the complexity of the problem from $\mathcal{O}(N^2)$ for a direct matrix multiplication implementation of the sum to $\mathcal{O}(N \log N)$. The Riemann sum approximation is useful when a discrete approximation to the continuous Fourier transform is required, and it may be preferable to the discrete Fourier transform in some cases because it is directly interpretable in the context of the Fourier transform. The Riemann sum approximation of the inverse Fourier transform applied to the Riemann sum approximation of the forward Fourier transform of a vector $\vec{x}$ returns the same vector $\vec{x}$—that is, the approximated transforms are still inverses of each other. However, the method pre-sented here results in the domain of the transformed function being defined only on integer multiples of $\frac{2\pi}{\Delta b N}$ where $\Delta$ is the spacing between values in the conjugate domain. Therefore, the transforms can only be implemented in such a way that they are each others' inverse if the input domain is defined on integer multiples of some number. Future work could try to alleviate this integer multiple constraint, or could try express higher-order integral approximations (like Simpson's rule) of the Fourier integral in terms of the discrete Fourier transform.

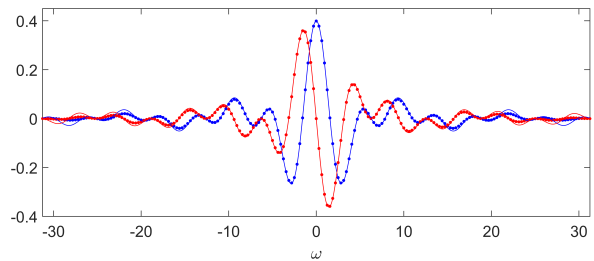

FIG. 2. The real (solid blue line) and imaginary (solid red line) parts of the analytic Fourier transform of $f(t)$ given by equation (9) are plotted with the real (blue ·) and imaginary (red ·) parts of the Riemann sum approximation as computed by MATLAB function `FT`. `FT` used 201 evenly-spaced samples of $f(t)$ on the interval $[-10, 10]$ with $a = 0$, $b = -1$. The approximation performs well for $|\omega| \ll \frac{1}{\tau'}$.

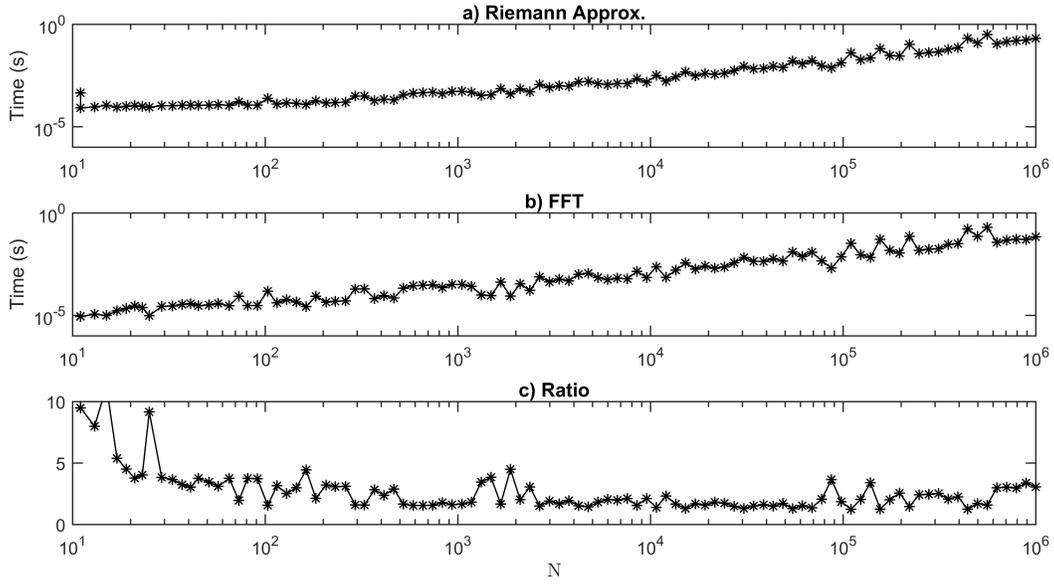

FIG. 3. Algorithm execution times as a function of input vector length for **a)** the Riemann approximation method as given in the function `FT` (see Section A1) and **b)** MATLAB's built-in fast Fourier transform function, `fft`. **c)** The ratio of times in (a) to (b). The ratio approaches a constant value as $N$ increases, indicating that the complexity of function `FT` scales similarly to `fft`—that is, as $\mathcal{O}(N \log N)$. The transforms were performed on the rectangular function $f(t) = \text{rect}(t-1)$, with sample points evenly-spaced on $[-10, 10]$.

## A1. Appendix

A simple implementation of the forward transform in MATLAB is given below. The function `FT` takes input vectors of length $N$ `t,f` and input scalars `w_1,a,b` (with `w_1` an integer multiple of $\frac{2\pi}{\tau' b N}$) and returns output vectors of length $N$ `w,ff` where `w(1) = w_1`.

```matlab
function [w,ff] = FT(t,f,w_1,a,b)

N = numel(t); % size of input vectors
tau = ((max(t)-min(t))./(N-1)); % input sampling period
W = -2*pi/(tau*b*N); % define frequency spacing
w_k = W.*(0:N-1); % define frequency vector beginning at zero

X = sqrt(abs(b)./((2*pi).^(1-a)))*tau.*exp(1i*b*w_k*t(1)).*fft(f);
% unshifted Riemann sum approximation using built-in "fft" function

w = w_k+w_1; % shift frequency vector so that w(1) = w_1

n = round(w_1/W); % number of indices that w is shifted from w_k by

m = zeros(size(w_k)); % initialize vector of periodic shift indices
if n > 0
    % determine m for each entry in w:
    idx = mod(n,N);
    m(N-idx:N) = floor(n./(N+1)) + 1;
    m(1:N-idx-1) = floor(n./(N+1));

    inds = circshift(1:N,-idx,2); % circularly shift indices of X to
    % corresponding new w positions
elseif n < 0
    % determine m for each entry in w:
    idx = mod(-n,N);
    m(1:idx) = -( floor(-n./(N+1)) + 1 );
    m(idx+1:N) = - floor(-n./(N+1));

    inds = circshift(1:N,idx,2); % circularly shift indices of X to
    % corresponding new w positions
end

phase = (exp(1i*(2*pi/tau)*t(1))).^m; % generate phases to shift by

ff = phase.*X(inds); % phase shift elements of circularly-shifted X, output
% answer
end
```

The corresponding inverse transform function `IFT` is also given for convenience:

```matlab
function [t,f] = IFT(w,ff,t_1,a,b)

N = numel(w); % size of input vectors
W = ((max(w)-min(w))./(N-1)); % input sampling period
tau = -2*pi/(W*b*N); % define time spacing
t_j = tau.*(0:N-1); % define time vector beginning at zero

x = sqrt(abs(b)./((2*pi).^(1+a)))*W.*exp(-1i*b*t_j*w(1)).*N.*ifft(ff);
% unshifted Riemann sum approximation using built-in "ifft" function

t = t_j+t_1; % shift time vector so that t(1) = t_1

n = round(t_1/tau); % number of indices that t is shifted from t_j by

m = zeros(size(t_j)); % initialize vector of periodic shift indices
if n > 0
    % determine m for each entry in w:
    idx = mod(n,N);
    m(N-idx:N) = floor(n./(N+1)) + 1;
    m(1:N-idx-1) = floor(n./(N+1));

    inds = circshift(1:N,-idx,2); % circularly shift indices of X to
    % corresponding new w positions
elseif n < 0
    % determine m for each entry in w:
    idx = mod(-n,N);
    m(1:idx) = -( floor(-n./(N+1)) + 1 );
    m(idx+1:N) = - floor(-n./(N+1));

    inds = circshift(1:N,idx,2); % circularly shift indices of X to
    % corresponding new w positions
end

phase = (exp(-1i*(2*pi/W)*w(1))).^m; % generate phases to shift by

f = phase.*x(inds); % phase shift elements of circularly-shifted x, output
% answer
end
```